# Walrasian Economic Equilibrium Problems in Convex Regions

**Renying Zeng**
Mathematics Department, Saskatchewan Polytechnic, Saskatoon, SK, Canada
renying.zeng@saskpolytech.ca

**Abstract.** This paper defines Walrasian economic equilibrium in convex regions, introduces the definition of strict-proper quasimonotone functions, and derives some necessary and sufficient conditions for the existence and uniqueness of Walrasian equilibrium vectors in convex regions.

## 1. Introduction

General Equilibrium Theory was developed by the French mathematical economist Leon Walras (1834-1910). Walras argues that, in a perfect competition economy, the participants in the market with a certain number of goods as a supply, will be able to achieve a balance.

He also believes that the total revenue of all production factors and the total revenue from the sale of all consumer goods will be equal under the "perfectly competitive" equilibrium. The essence of this theory is that the capitalist economy can be in a stable equilibrium state. In the capitalist economy, consumers can get the most out of it, the entrepreneur can get the maximum profit, and the owner of the production factors can get the maximum reward.

Walras argues that the process of achieving a balanced price or value is consistent, so price decisions and value decisions are one thing. He believes that various goods, services supplies, demand quantities, and prices are interrelated. A commodity price and quantity changes can cause other goods quantity and price changes. So we can not only study a commodity, a market supply and demand changes, we must also study all the goods, all the market supply and demand changes. Only when all markets are in balance, individual markets can be in a state of equilibrium.

Walrasian general equilibrium price decision thought is illustrated by “mathematical formulas”. Through the solutions of equations, Walras proved that there were a series of market prices and the number of transactions in the market (these prices and quantities are the equilibrium prices and quantities), so that each consumer, entrepreneur and resource owner achieve their purpose, so that society can exist harmoniously and steadily.

This paper works with Walrasian equilibrium problems as mathematical formulas—variational inequalities.

## 2. Preliminaries

Let $X$ be a non-empty subset of the n-dimensional Euclidean space $R^n$ and let $E$: $X \to R^n$ be a function. A point $x^* \in X$ is said to be a solution of the **variational inequality** $VI$ ($E$, $X$ ) if there holds

$$E(x^*)^T(x - x^*) \geq 0, \forall x \in X. \tag{1}$$

A function $E$: $X \to R^n$ is said to be **pseudomonotone** if there holds

$$E(x)^T(y - x) \geq 0 \Rightarrow E(y)^T(y - x) \geq 0, \forall y \in X;$$

**strictly pseudomonotone** if there holds

$$E(x)^T(y - x) \geq 0 \Rightarrow E(y)^T(y - x) > 0, \forall y \in X,\ y \neq x.$$

Similar to the proof in [1], one has the following Lemma 1.

**Lemma 1** Let $X \subseteq R^n$ be a convex region and $E$: $X \to R^n$ a continuous function, then the variational inequality $VI[E, X]$ has solutions.

**Lemma 2** Let $X \subseteq R^n$ be a convex region and $E$: $X \to R^n$ a continuous and strictly pseudomonotone function, then the variational inequality $VI[E, X]$ has a unique solution.

## 3. An Example of Walrasian Economic Equilibrium Problem

This section presents an example of a Walrasian Equilibrium problem in a convex region.

**Walrasian equilibrium.** Consider a perfect competition economy with $n$ commodities. Give a price vector $p \in R^n$, the **aggregate excess demand function** is defined by

$$\begin{aligned} E(p) &= -D(p) + S(p) \\ &= (E_1(p), E_2(p), \ldots, E_n(p)), \end{aligned}$$

where $D$ and $S$ are the **demand** and **supply functions**, respectively.

Traditionally, a vector $p^*$ is said to be a **Walrasian equilibrium price vector** [1], [2] if it solves the following variational inequality

$$E(p^*) \geq 0, \forall p \in \mathsf{S}.$$

Where

$$\mathsf{S} = \{x \in \mathsf{R}^n; x_j \geq 0, \sum_{j=1}^{n} x_j = 1\}$$

is an $(n-1)$-dimensional **unit simplex**.

It has been known [3] Theorem 7.2 that a vector $p^*$ is a Walrasian equilibrium price vector if and only if it is a solution of the variational inequality

$$E(p^*)^T (p - p^*) \geq 0, \forall p \in \mathsf{S}.$$

If each price of a commodity is involved in the market structure which has a lower positive bound and may have an upper bound, then the prices are assumed to be contained in the **box-constrained set**

$$\mathsf{K} = \{x \in \mathsf{R}^n; 0 \leq a_j \leq x_j \leq b_j \leq +\infty\},$$

where $a_j$ are constants, $b_j$ are either constants or $+\infty$ ($j = 1, \ldots, n$). Then, a vector $p^*$ is a Walrasian equilibrium price vector if it is a solution of the variational inequality over a box-constrained set $\mathsf{K}$ [3]:

$$E(p^*)^T (p - p^*) \geq 0, \forall p \in \mathsf{K}.$$

Now, the following Walrasian equilibrium problem over a convex feasible region is introduced, which extends the above definition that over a box-constrained set.

**Definition 1** A price vector $p^*$ is a Walrasian equilibrium price vector if it is a solution to the variational inequality

$$E(p^*)^T (p - p^*) \geq 0, \forall p \in \mathsf{X}. \tag{2}$$

Where $\mathsf{X}$ is a convex region.

An example is given here to illustrate that an economic equilibrium problem may have a convex feasible region.

**Example 1** The function

$$A = P \frac{p_1 (1 + p_1)^N}{(1 + p_1)^N - 1}$$

is known as the mortgage payment amount that should be paid periodically for $N$ periods on a mortgage amount $P$ at a periodic interest rate of $p_1$. The mortgage can be considered as a type of commodity. The price of the mortgage would be the interest rate $p_1$, which means $p_1$ dollars for each one hundred-mortgage for one period. At the right time of the mortgage has been paid off, the value $p_2$ of the house satisfies

$$p_2 \geq B + A = B + P\frac{p_1(1+p_1)^N}{(1+p_1)^N - 1},$$

where $B$ is the down payment in the beginning of the mortgage.

Consider the function

$$p_2 = f(p_1) = B + P\frac{p_1(1+p_1)^N}{(1+p_1)^N - 1}.$$

If $N = 2$, then

$$p_2 = B + P\frac{p_1(1+p_1)^2}{(1+p_1)^2 - 1} = B + P\frac{p_1^2 + 4p_1 + 3}{p_1 + 2}.$$

Then the derivative

$$p_2' = f'(p_1) = P\frac{p_1^2 + 4p_1 + 3}{(p_1 + 2)^2} > 0, \forall p_1 \geq 0$$

So, the function $p_2 = f(p_1)$ is increasing if $p_1 > 0$.

Then,

$$p_2'' = f''(p_1) = P\frac{2}{(2 + p_1)^3} > 0, \forall p_1 \geq 0$$

deduces that $p_2 = f(p_1)$ is concave up if $p_1 \geq 0$.

Then, the inverse function $p_1 = g(p_2)$ of $p_2 = f(p_1)$ is increasing and concave up, since the graphs of $p_1 = g(p_2)$ of $p_2 = f(p_1)$ are symmetric to the line $p_2 = p_1$.

Similar discussions can be applied to the cases of $N > 2$, for the interest rate $p_1$ varies in a certain range.

For a perfect competition economy with $n$ commodities, let

$$X = \{p \in R^n; 0 \le p_1 \le g(p_2),\ 0 \le a_j \le p_j \le b_j \le +\infty,\ j = 1,2,...n\}$$

where $p = (p_1, p_2, ..., p_n) \in R^n$ is the price vector. If $b_j < +\infty$ ($j =$ 2, .., $n$), then X is a convex region.

Now we know that some Walrasian economic equilibrium problems in a perfect competition economy can be formulated as variational inequalities over convex regions.

## 4. Main Results

A price vector $p^* \in X$ is said to be a **stationary price vector** of the Walrasian equilibrium problem (2) if there holds

$$E(p)^T (p - p^*) \ge 0, \forall p \in X. \quad (3)$$

A. Daniilidis and N. Hadjisavvas [9] introduced the following definition of the proper quasimonotonicity.

The function $E$: $X \to R^n$ is said to be **properly quasimonotone** if $\forall p_i \in X$, $j \in \{1, 2, \ldots, m\}$, $\forall p \in con(p_1, p_2, \ldots, p_m)$, $\exists j \in \{1, 2, \ldots, m\}$ such that

$$E(p_j)^T (p_j - p) \ge 0.$$

The following Theorem 1 gives a necessary and sufficient condition for the existence of the Walrasian equilibrium price vector.

**Theorem 1** There exists a stationary price vector $p^*$ for the Walrasian equilibrium problem (2) if and only if $E$: $X \to R^n$ is properly quasimonotone.

**Proof**. The proof of the sufficient is from A. Daniilidis and N. Hadjisavvas [10], and the necessity can be found in R. John [11].

We introduce the definition of **strict-proper quasimonotonicity**.

**Definition 2** The function $E$: $X \to R^n$ is said to be **strict-properly quasimonotone** if $\forall p_i \in X$, $i \in \{1, 2, \ldots, m\}$, $\forall p \in con(p_1, p_2, \ldots, p_m)$ with $p \ne p_i$, $i \in \{1, 2, \ldots, m\}$, $\exists j \in \{1, 2, \ldots, m\}$ such that

$$E(p_j)^T (p_j - p) > 0.$$

**Theorem 2** If $E$: $X \to R^n$ is strict-properly quasimonotone, then $E$ is strictly pseudomonotone, i.e., $\forall p, q \in X$, $p \ne q$,

The following Theorem 3 gives necessary and sufficient conditions for the existence and uniqueness of the Walrasian equilibrium price vector and of the stationary price vector.

**Theorem 3** Let $E$: $X \to R^n$ be continuous. Then the following are equivalent:

(a) $E$: $X \to R^n$ is strict-properly quasimonotone;
(b) There exists a unique stationary price vector $p^*$ for the Walrasian equilibrium problem (2);
(c) There exists a unique Walrasian equilibrium price vector for the Walrasian equilibrium problem (2).

## 5. Conclusion

This paper obtains some necessary and sufficient conditions for the existence and uniqueness of the solution of a Walrasian equilibrium problem over a convex region. Which means that, in a pure exchange economy (perfect competition economy), the participants in the market with a certain number of goods as a supply, will be able to achieve a balance.